\documentclass{article}

\usepackage{amssymb,amsmath,latexsym,amsthm}
\usepackage{amsgen}
\usepackage{amsfonts}

\newtheorem{theorem}{Theorem}[section]
\newtheorem{lemma}[theorem]{Lemma}

\newtheorem{proposition}[theorem]{Proposition}

\theoremstyle{definition}
\newtheorem*{definition}{Definition}
\newtheorem*{remark}{Remark}

\newtheorem{algorithm}{Algorithm}

\def\qqq{\mathbb{Q}}
\def\rrr{\mathbb{R}}

\def\zzz{\mathbb{Z}}

\def\pf{\textbf{ proof}:\ }
\def\qed{$\Box$}

\begin{document}

\title{Continued fractions and Parallel SQUFOF}


\author{S. McMath,
F. Crabbe\thanks{
U. S. Naval Academy, Computer Science Department,
Annapolis, MD 21402, crabbe@usna.edu},
D. Joyner\thanks{
U. S. Naval Academy,
Mathematics Department,
Annapolis, MD 21402, wdj@usna.edu}}

\date{1-9-2006}
\maketitle


\small
        \vspace*{.1in}
        \begin{center}%
          {\bfseries Dedication\vspace{-.5em}\vspace{.1in}}%
	\ \ To the memory of Daniel Shanks (1917-1996).
        \end{center}%

\begin{abstract}
In this partly expository paper, we prove two results.
\begin{itemize}
\item
That the two-sided continued fraction
of the normalized square root (an important part of the SQUFOF
algorithm) has several very attractive properties - periodicity,
a symmetry point corresponding to a factorization of $N$, and so on.
\item 
The infrastructure distance formula.
\end{itemize}
Finally, we describe a method for parallelization of SQUFOF 
that maintains its efficiency per procesor as the number of 
processors increases, and thus is predicted to be useful 
for very large numbers of processors.
\end{abstract}

\section{Introduction}

Though there are many fast algorithms
for factoring numbers, this paper focuses on one known as square forms
factorization or \emph{SQUFOF} (see Algorithm \ref{alg:fact} below for
a precise description).  Daniel Shanks developed SQUFOF in the 1970's,
and it is still the fastest known algorithm for factoring integers in
the $20$- to $30$-digit range.  SQUFOF is used to this day in
conjunction with other factorization algorithms that need to factor
$20$-digit numbers in order to generate the factors of higher digit
numbers.

Most of the Shanks' original work on SQUFOF was not published (see
however, \cite{Sh1})
and his notes are incomplete\footnote{These notes have been typed in
LaTeX and are available on the web \cite{Sh2}, \cite{Sh3},
\cite{Sh4}.}.  One purpose of this paper
is to present Shanks's original SQUFOF algorithm in its entirety.  
The paper goes on to present several results concerning 
both traditional SQUFOF and its parallelization.

This paper contains three main results:
\begin{enumerate}

\item 
A proof that the two-sided continued fraction
of the normalized square root (an important part of the SQUFOF
algorithm) has several very attractive properties - periodicity,
a symmetry point corresponding to a factorization of $N$, and so on
(see Theorems \ref{RieselProof}, \ref{strongPeriod}, and
\ref{Symmetries} for details).
This result was probably known 200 years ago to Lagrange
and Galois and Gauss - see for example, Perron \cite{P},
Buell \cite{B}, and Williams \cite{W}.

\item 
A proof of the infrastructure distance formula, 
Theorem \ref{thrm:14} below, which is
also an important part of SQUFOF. This is in some sense well-known
but a proof has not, as far as we can see, appeared in the
literature.
(However, see Cohen \cite{Coh2}, Proposition 5.8.4, and 
Williams and Wunderlich \cite{WW} Theorem 5.2 for closely
related results.)

\item 
Investigation of a method for parallelization of 
SQUFOF that maintains its efficiency
per procesor as the number of processors increases, and thus is
predicted to be useful for large numbers of processors. (See
\cite{W}, \cite{WW}, and \cite{Gower} for work on similar ideas.)
The implementation in C, and subsequent numerical data 
due to the first author, is new as far as we know. This is
briefly sketched in \S 4 below.

\end{enumerate}

\vskip .1in
Although the theoretical results in this paper are known to the experts, it
is hoped that putting all these results together will serve a
useful purpose. 
This paper is a version of the first author's 
undergraduate ``Trident'' thesis, advised by the second two authors.

\section{Continued Fractions and Quadratic Forms}

The stepping stone for SQUFOF is the continued fraction expansion for
the square root of $N$. (We slightly simplify matters by instead using
the ``normalized square root (equation \ref{eqn:normsqrt}) here.)  The
terms of this continued fraction expansion give rise to a sequence of
quadratic forms of discriminant $N$ via (\ref{eqn:CF2QF}).  We shall
describe SQUFOF in terms of the ``cycle'' of continued fractions in
the periodic expansion of (\ref{eqn:normsqrt}) and the corresponding
quadratic forms.

\subsection{Integral binary quadratic forms}

There is a ``dictionary'' between certain aspects of
\begin{itemize}
\item indefinite integral binary quadratic forms,
\item ideals in a real quadratic number field,
\item the simple continued fraction of quadratic surds.
\end{itemize}
The reader will be assumed to be familiar with at least the basic
aspects of this correspondence.  For details, see for example, Buell
\cite{B}, Lenstra \cite{Len}, Williams \cite{W} (especially
pp. 641-645), Cohen \cite{Coh} and the references found there,
or \cite{M}.

A \textbf{binary quadratic form} (or simply a ``form'') is a homogeneous
form of degree two in two variables $x,y$,

\[
f(x,y)=ax^2+bxy+cy^2=
( x,y) \cdot
\left(
\begin{array}{cc}
a & b/2\\
b/2 & c\end{array}
\right)\cdot 
\left(
\begin{array}{c}
x\\
y\end{array}
\right),
\]
for some constants $a,b,c$. This form shall also be denoted
by the triple $(a,b,c)$. The \textbf{discriminant}\footnote{Sometimes 
also called the \textbf{determinant} of $f$.}
of $f$ is $D=disc(f)=b^2-4ac$.
We shall focus on the case $D>0$, 
in which case the form
is called \textbf{indefinite}. 
{\it From now on, we assume without further
mention that $D>0$ is a non-square such that 
$D \equiv 0\pmod 4$ or $ D \equiv 1\pmod 4$.}

If $a,b,c \in \zzz$ then we say $f$ is \textbf{integral}.  If moreover
$gcd(a,b,c)=1$, then we say the form is \textbf{primitive}.  Let
$F(D)$ denote the set of all integral forms of discriminant $D$ and
let $F(D)_{p}$ denote the subset of primitive ones.

The groups

\[
GL_2(\zzz)=
\{
\gamma=\left(
\begin{array}{cc}
s & t\\
u & v\end{array}
\right)\ |\ s,t,u,v\in \zzz,\ \det(\gamma)=\pm 1\},
\]
and
\[SL_2(\zzz)=
\{
\gamma\in GL_2(\zzz)\ |\ \det(\gamma)=1\}
\]
act on the polynomials $\zzz[x,y]$ via

\[
\gamma =\left(
\begin{array}{cc}
s & t\\
u & v\end{array}
\right): (x,y)\longmapsto 
(sx+ty,ux+vy).
\]
Therefore, they also act on the set of integral forms 
via

\[
(\gamma^* f)(x,y)=f(sx+ty,ux+vy),
\]
for $\gamma \in GL_2(\zzz)$.
In terms of the symmetric matrix
$A=\left(
\begin{array}{cc}
a & b/2\\
b/2 & c\end{array}
\right)$
associated to the form $f$,
this action may be epressed as

\[
\gamma^*(A)=\, ^t\gamma\cdot A\cdot \gamma.
\]
We say that two forms $f_1,f_2$ are \textbf{equivalent}
if $f_2=\gamma^*f_1$, for some $\gamma\in GL_2(\zzz)$.
We say that two forms $f_1,f_2$ are \textbf{properly equivalent},
written $f_1\sim f_2$,
if $f_2=\gamma^*f_1$, for some $\gamma\in SL_2(\zzz)$.
For $f\in F(D)$, we let

\[
F(D)_{f}=[f]=\{f'\in F(D)\ |\ f\sim f'\}
\]
denote the proper equivalence class of $f$.
An element $\gamma\in GL_2(\zzz)$ is called an \textbf{automorph} of
$f$ if $\gamma^*f=f$. A form $f$ is called \textbf{ambiguous} if 
it has an automorph in $GL_2(\zzz)-SL_2(\zzz)$. 
Note that if $f\in F(D)$ is ambiguous then 
each $f'\in [f]$ is also ambiguous. 

We say that two forms
$(a_1,b_1,c_1),(a_2,b_2,c_2)\in F(D)$ are 
\textbf{adjacent} if $c_1=a_2$ and
$b_1+b_2\equiv 0 \pmod{2a_2}$. In this case, we say that 
$(a_2,b_2,c_2)$ is to the \textbf{right of} 
$(a_1,b_1,c_1)$ ($(a_1,b_1,c_1)$ is to the \textbf{left of}
$(a_2,b_2,c_2)$).

\subsubsection{Reduction}

A form $(a,b,c)$ is called \textbf{reduced}
if $|D^{1/2}-2|a||<b<D^{1/2}$.
Let $F(D)_{r}$ denote the subset of
reduced forms of discriminant $D$.

\begin{lemma}\label{buell-lemma}
(a) Given any $f\in F(D)_{r}$ there is a unique
$f'\in F(D)_{r}$ adjacent to the right of $f$
and a unique $f''\in F(D)_{r}$ adjacent to the left of $f$.

(b) There are exactly two 
reduced ambiguous forms in a cycle of reduced forms in an
ambiguous class.
\end{lemma}

For (a) see Buell \cite{B}, page 23; for (b), see \cite{B}, Theorem
9.12.  Lemma \ref{buell-lemma} allows us to define the \textbf{cycle} of
reduced forms associated to $f\in F(D)_{r}$: it is the set of all
$f'\in F(D)_{r}$ which is adjacent to the left or right of $f$.
This cycle is denoted $F(D)_{r,f}$.

\begin{lemma}
\label{2Symmetries}
An ambiguous equivalence class contains two points of symmetry, that
is, pairs of reduced adjacent forms, $(c,b,a)$ to the left of
$(a,b,c)$, in the cycle that are the symmetric reverse of each other.
In that case, either $a$ divides the determinant, or $a/2$ divides the
determinant.
\end{lemma}

This follows from Theorem \ref{Symmetries} below.

It is evident that if a form is ambiguous, then each form
in its equivalence class is also ambiguous.  

\begin{proposition}
\label{prop:cycles}
The set $F(D)_{r}$ of reduced forms can be
partitioned into cycles of adjacent forms.
\end{proposition}

Consider the action of 

\[
T_m=\left(
\begin{array}{cc}
1 & m\\
0 & 1\end{array}
\right)
\]
on a form $(a,b,c)$:
$T_m(a,b,c)=(a',b',c')$,
where $a'=a$, $b'=b+2am$, $c'=\frac{(b')^2-D}{4a'}$.
This defines a map
$T_m:F(D)\rightarrow F(D)$, for each $m\in \zzz$.

Consider the action of 

\[
W=\left(
\begin{array}{cc}
0 & -1\\
1 & 0\end{array}
\right)
\]
on a form $(a,b,c)$:
$W(a,b,c)=(a',b',c')$,
where $a'=c$, $b'=-b$, $c'=a$.
This defines a map $W:F(D)\rightarrow F(D)$.

\begin{algorithm}
\label{alg:red}
{\small{
{\fontfamily{phv}\selectfont
(Reduction)

Input: $f\in F(D)$.

Output: $f'\in F(D)_{r}$ with $f\sim f'$.

Let $f(x,y)=ax^2+bxy+cy^2$ and let

\[
J_{a,D}=\{x\ |\ -|a|<x<|a|,\ {\rm if}\ |a|\geq D^{1/2},\ 
-2|a|<x<D^{1/2}, \ {\rm if}\  |a|<D^{1/2}\}.
\]

\begin{enumerate}
\item
Apply $T_m$ to $(a,b,c)$ to obtain a form
$(a,b',c')$, where $b'\in J_{a,D}$ and $c'$ is chosen so that
the new form has discriminant
$D$.
\item \vspace{-3mm}
If $(a,b',c')$ is reduced then return $f'(x,y)=ax^2+b'xy+c'y^2$.
Otherwise, replace $(a,b',c')$ by 
$W(a,b',c')=(c',-b',a)$ and go to step 1.

\end{enumerate}
}}}
\end{algorithm}
According to Lagarias \cite{L1}, this has complexity
$O(\log (max(|a|,|b|,|c|)))$.

Define the \textbf{adjacency map} $\rho :F(D)\rightarrow F(D)$
by

\begin{equation}
\label{eqn:adj}
\rho(a,b,c)=(a',b',c'),
\end{equation}
where $a'=c$, $b'\in J_{c,D}$, and
$b'\equiv -b \pmod{2c}$, and $c'$ is determined by the condition
$disc(a',b',c')=D$. This defines a bijection
$\rho:F(D)_{r}\rightarrow F(D)_{r}$.

Unfortunately, given $f\in F(D)$ with $D>0$ there are usually several
$f'\in F(D)_{r}$ which are properly equivalent to $f$.
In other words, the cycle

\[
F(D)_{r,f}=\{f'\in F(D)_{r}\ |\ f\sim f'\}
=\{f'=\rho^nf \ |\ n\in \zzz\}
\]
can be rather large. Indeed, it is known that
$
|F(D)_{r,f}|=O (D^{1/2+\epsilon}), 
$ (where the $O$-constant depends on $\epsilon$)
for all $\epsilon >0$, where the exponent $1/2$ is best possible
(Lagarias \cite{Len,L2}) and .

\subsubsection{Composition}

The \emph{composition} of forms has important properties for SQUFOF.
The rules of composition are fairly general. A binary quadratic form
$F$ is called \textbf{a composition} of $f,g\in F(D)$ if it satisfies an
equation such as

\begin{equation}
\label{eqn:comp}
f(x,y)g(u,v)=F(B_1(x,y,u,v),B_2(x,y,u,v)),
\end{equation}
where $B_1$ and $B_2$ are quadratic forms in
$x,y,u,v$ of a certain type. The exact conditions $B_1,B_2$ satisfy
do not concern us here (see Cox \cite{Cox} if
you are curious and Gauss \cite{Gauss} if you are really
curious). The point is that there may be more than one pair
$B_1,B_2$ satisfying (\ref{eqn:comp}), so that the composition
$F$ is not unique. (However, the conditions
on $B_1,B_2$ specified by Gauss do imply that, for a given
$f,g\in F(D)$ any two such compositions must be equivalent to 
each other.) One way around this ambiguity is to specify a choice 
of $B_1, B_2$ and hence define $F$ uniquely. 

%

The idea described below was known in some form to Dirichlet 
and possibly Gauss.

\begin{algorithm}
\label{alg:comp}
\noindent
{\small{
{\fontfamily{phv}\selectfont
Input: $(a_1,b_1,c_1), (a_2,b_2,c_2)\in F(D)$.

Output: A composition 
$(\frac{a_1a_2}{m^2},B,\frac{(B^2-D)m^2}{4a_1a_2})\in F(D)$.

\begin{enumerate}
\item
Compute 
$m=gcd(a_1,a_2,\frac{b_1+b_2}{2})$.
(Since $D=b_i^2-4a_ic_i$, for $i=1,2$,
$b_1$ and $b_2$ have the same parity.)

\item \vspace{-3mm}
Solve the congruences

\[
\begin{array}{c}
a_2mB\equiv mb_1a_2 \pmod{2a_1a_2} ,\\
a_1mB\equiv mb_2a_1 \pmod{2a_1a_2} ,\\
\frac{b_1+b_2}{2}mB\equiv m\frac{b_1b_2+D}{2} \pmod{2a_1a_2} ,\\
\end{array}
\]
simultaneously an integer $B$. Choose the solution with smallest
absolute value.
\end{enumerate}

}}}
\end{algorithm}

See \cite{Sh1} or \cite{B} for a proof of the correctness of this 
algorithm. Buell \cite{B} also provides the substitutions that would be 
needed for Gauss's definition of composition. 

In other words, we define the \textbf{composition} of 
$(a_1,b_1,c_1), (a_2,b_2,c_2)\in F(D)$
to be the form resulting from the above algorithm:

\[
(a_1,b_1,c_1)* (a_2,b_2,c_2) 
=(\frac{a_1a_2}{m^2},B,\frac{(B^2-D)m^2}{4a_1a_2}).
\]

\begin{remark}
\label{remark:comp}
The binary operation $*:F(D)\times F(D)\rightarrow F(D)$ is 
associative but not its ``restriction''
$\#:F(D)_{r}\times F(D)_{r}\rightarrow F(D)_{r}$
(where $\#$ is composition algorithm \ref{alg:comp}
followed by reduction algorithm \ref{alg:red}).

Let $f,g\in F(D)_{r}$ be elements in
the principal cycle of discriminant $D$. 
It was observed by Shanks (see \S 5 in Lenstra \cite{Len})
that cycles enjoy a ``coset-like property''
$\rho^k f\#\rho^\ell g=\rho^{a_{k,\ell}}(f\# g)$, for some
$a_{k,\ell}\in \zzz$ . In particular, the principal cycle is
closed under composition. Therefore, the 
the set of complete quotients of the continued fraction of 
such an $\alpha$ can be identified with a set closed under
$\#$.

For further discussion of this, see Lenstra \cite{Len} (5.1).

The ``structure'' of a cycle has been termed 
the ``infrastructure'' of $F(D)$ by Shanks.

\end{remark}

If $f,f',g,g',h\in F(D)$ then Gauss showed

(a) $(f*g)*h\sim f*(g*h)$, and 

(b) $f\sim f'$ and $g\sim g'$ imply $f*g\sim f'*g'$.

\noindent
These imply that the set of equivalence classes of forms of
discriminant $D$ is a group $C(D)$, called the \textbf{class group} of
$D$. From the construction, it is clear that $f*g\sim g*f$, so $C(D)$
is abelian.

The following Theorem was known to Shanks, since SQUFOF depends 
essentially on it.

\begin{theorem}
An equivalence class has order $2$ or $1$ in the class group if and only 
if it is ambiguous.
\label{Order2}
\end{theorem}

Any form $(1,b,c)\in F(D)$ acts as the identity for $*$.
The cycle of the identity is the \textbf{principal cycle} of forms.
Any form $f$ whose square $f^2=f*f$ belongs to the
principal cycle is an ambiguous form
(\cite{B}, Corollary 4.9).

\subsection{Continued fractions}
\label{Continued Fractions}

Throughout, assume that $N\equiv 1 \pmod 4$ and is 
not a perfect square.  

We shall only consider simple continued fractions here.
In other words, if $\alpha\in \rrr$ is the number we want 
to compute the continued fraction of, let $x_0 = \alpha$, 
$b_0 = \left\lfloor \alpha \right\rfloor$,
where $\left\lfloor x\right\rfloor$ denotes the \textbf{floor} of $x$,
and, for $i>0$, let

\begin{equation}
\label{expand}
x_i = \frac{1}{x_{i-1}-b_{i-1}},\ \ \ \ b_i = \left\lfloor x_i \right\rfloor .
\end{equation}
The term $x_i$ is called the \textbf{$i^{th}$ complete quotient} of
$\alpha$ and $b_i$ is called the \textbf{$i^{th}$ partial quotient} of
$\alpha$.
The \textbf{simple continued fraction} of
$\alpha$ is (\cite{H&W}):
\[
\alpha = b_0 + \frac{1}{b_1+\frac{1}{b_2+...}} ,
\]
also written $[b_0,b_1,b_2,...]$.
We are only concerned with continued fractions of an
irrational $\alpha \in K=\qqq (\sqrt{N})$. In this case, the sequence
$b_0,b_1,b_2,...$ is eventually periodic. 

For example, let 

\begin{equation}
\label{eqn:normsqrt}
\alpha = 
\left\{
\begin{array}{cc}
\frac{\sqrt{N}+\left\lfloor \sqrt{N} \right\rfloor-1}{2}, 
& \left\lfloor \sqrt{N} \right\rfloor \ {\rm even},\\
\frac{\sqrt{N}+\left\lfloor \sqrt{N} \right\rfloor}{2}, 
& \left\lfloor \sqrt{N} \right\rfloor \ {\rm odd}.
\end{array}
\right.
\end{equation}
We call this $\alpha$ the \textbf{normalized square root of $N$}.
The continued fraction sequence $b_0,b_1,...$ is (purely) periodic.
In general, the period of $\alpha$ is the size of the cycle
associated to the identity in the class group (Buell \cite{B},
Theorem 3.18 (a)).

At each step in the continued fraction expansion, 
it is possible to simplify $x_i-b_i$ to the form 
$\frac{\sqrt{N}-P_i}{Q_i}\in [0,1)$, 
where $P_i,Q_i\in \zzz$ satisfy $P_i^2 \equiv N \pmod {Q_i}$.
In general, if $P,Q$ are positive integers and
$x=\frac{\sqrt{N}+P}{Q}$ satisfies
$P^2\equiv N\pmod{Q}$, $0<P<\sqrt{N}$,
$|\sqrt{N}-Q|<P$, then we say that $x$ is
\textbf{reduced}. It is known that if
$x,y$ are two such reduced numbers and $y=\gamma (x)$
(where $\gamma = \left(\begin{array}{cc}
a & b\\
c & d\end{array}\right) \in  SL_2(\zzz)$ acts on $\hat{\rrr} =
\rrr \cup \{\infty\}$ by $\gamma (x) = \frac{a x + b}{c x + d}$)
then $y$ occurs in the simple continued fraction 
expansion of $x$ as a complete quotient
(and $x$ occurs in the simple continued fraction 
expansion of $y$ as a complete quotient). See Buell \cite{B},
Proposition 3.20 for a proof.

If $P,Q$ are positive integers and
$x=\frac{\sqrt{N}+P}{Q}$ then we associate to $x$ the
quadratic forms

\begin{equation}
\label{eqn:CF2QF}
f_-=(-Q/2,P,-\frac{P^2-N}{2Q}),\ \ \ 
f_+=(Q/2,P,\frac{P^2-N}{2Q}),
\end{equation}
which have discriminant $N$. (We implicitly assume here that 
$\frac{P^2-N}{2Q}\in \zzz$ and $Q$ is even.
Note that if $x$ is reduced then so are
$f_\pm$, and conversely.)

\begin{lemma} (H. Cohen \cite{Coh}, \S 5.7.1)
The continued fraction expansion of the
quadratic irrational corresponding to the unit
reduced form is not only periodic but symmetric.
\end{lemma}

What is the continued fraction analog of ``adjacency'' of forms?
Applying the adjacency map (\ref{eqn:adj}) is roughly analogous to the
``stepping'' process of going from one complete quotient to the next
in a continued fraction. See Williams \S 5 for a discussion of the the
ideal-theoretic analog, at least for the case of the simple continued
fraction of $\frac{-1+\sqrt{N}}{2}$.

One tool used by many different algorithms is the continued 
fraction expression for (\ref{eqn:normsqrt}), 
where $N$ is the number to be factored.  
This expression is calculated recursively:
$x_0 = \alpha$, $b_0 = \left\lfloor x_0 \right\rfloor$,
and using (\ref{expand}) in general.
Observe that solving equation (\ref{expand}) for $x_{i-1}$ gives 
$x_{i-1} = b_{i-1}+\frac{1}{x_i}$.  
  
The recursive formulas are, for $i \geq 0$,

\begin{equation}
\label{basicStep}
\begin{array}{ll}
x_{i+1} &= \frac{1}{x_i-b_i}\\
&= \frac{Q_i}{\sqrt{N}-P_i} \\
&= \frac{\sqrt{N}+P_i}{Q_{i+1}} \\
&= b_{i+1} + \frac{\sqrt{N}-P_{i+1}}{Q_{i+1}}, \\
b_i &= \lfloor x_i \rfloor.
\end{array}
\end{equation}
Theorem \ref{RieselProof} provides some well-known fundamental 
properties and identities of continued fractions.  

\begin{theorem} (\cite{Riesel})

In the continued fraction expansion of (\ref{eqn:normsqrt}), 
with $x_0=\alpha$, each $x_i$ reduces to the form 
$\frac{\sqrt{N}+P_{i-1}}{Q_i}$, with unique
$Q_i,P_i\in \mathbb{Z}$ satisfying

\begin{itemize}
\item[(a)]
 $N=P_i^2+Q_iQ_{i+1}$, 

\item[(b)]
 $P_i = b_iQ_i-P_{i-1}$, 

\item[(c)] 
$b_i = \left \lfloor \frac{\lfloor \sqrt{N} \rfloor 
+ P_{i-1}}{Q_i} \right \rfloor \geq 1$, 

\item[(d)] $0<P_i<\sqrt{N}$, 

\item[(e)] $|\sqrt{N}-Q_i| < P_{i-1}$, 

\item[(f)] $Q_i$ is an integer, 

\item[(g)] $Q_{i+1} = Q_{i-1} + b_i(P_{i-1}-P_i)$.  

\item[(h)] 
This sequence is eventually periodic.

\item[(i)]
$
\left \lfloor \frac{\sqrt{N}+P_i}{Q_i} \right \rfloor 
= \left \lfloor \frac{\sqrt{N}+P_{i-1}}{Q_i} \right \rfloor = b_i.
$
\end{itemize}
\label{RieselProof}
\end{theorem}

These denominators $\{Q_i\}$ will be referred 
to as \textbf{pseudo-squares}.  
(Indeed, for $i \geq 0$, if we write $[b_0,b_1,...b_i] = \frac{A_i}{B_i}$ 
then $A_{i-1}^2-B_{i-1}^2N = (-1)^iQ_i$
and so $A_{i-1}^2 \equiv (-1)^iQ_i \pmod N$.) 

\begin{remark}
The fact that each $x_i$ reduces to the form 
$\frac{\sqrt{N}+P_{i-1}}{Q_i}$ is important for computational 
efficiency because this together with (c) imply that floating 
point arithmetic is not necessary for any of these calculations.  
Also, by use of (b) and (g), the arithmetic used in this 
recursion is on integers $< 2\sqrt{N}$.
\end{remark}

Since the continued fraction is eventually periodic, it is reasonable
to consider that when it loops around on itself, the terms being
considered may have come from some terms ``earlier" in the recursion.
Lemma \ref{reversal} shows that by exchanging these two related
expressions, the direction is reversed.  The algorithm for stepping a
continued fraction expansion in the opposite direction will be
precisely the same as the one for the forward direction, except that
the numerator is changed first.  Note that this same change (with the
exception of $c_0$) could be achieved by merely changing the sign of
$P_{i-1}$.

\begin{lemma}
Let $N$, and, for $i\geq 0$, let $x_i,b_i,P_i,Q_i$
be as in Theorem \ref{RieselProof}.  
Let $y_0 = \frac{\sqrt{N}+P_{i+1}}{Q_{i+1}}$ and let 
$c_0 = \lfloor y_0 \rfloor$.  
If we define, for $j\geq 1$, $y_j 
= \frac{1}{y_{j-1}-c_{j-1}}$, $c_{j-1}=[y_{j-1}]$
then $c_0 = b_{i+1}$ and $y_j 
= \frac{\sqrt{N}+P_{i-j+1}}{Q_{i-j+1}}$, when
$0 \leq j\leq i$.
\label{reversal}
\end{lemma}

Using Lemma \ref{reversal} to go backwards in the
continued fraction expansion, denote the terms 
before $x_0$ as $x_{-1},x_{-2},...$.
The sequence $\{x_i\ |\ i\in \zzz\}$ will be called
the \textbf{two-sided continued fraction} of $x_0$.  
Define $Q_{-i}$ and $P_{-i}$ similarly, $i\geq 0$. 

\begin{theorem}
\label{strongPeriod}
(a) With these conventions on the negative indices,
Theorem \ref{RieselProof} applies for all $i\in \zzz$.

(b) Define $x_i$ as in Theorem \ref{RieselProof}, $i\in \zzz$.
There exists a positive integer $\pi$ such that for all $i\in \zzz$,
$x_i = x_{i+\pi}$.

(c) Let $x_0 = \alpha$ such that $Q_0 \mid 2P_{-1}$ (as in equation
(\ref{eqn:normsqrt})).  The sequence of pseudo-squares is symmetric
about $Q_0$, so that for all $i\in \zzz$, $Q_i = Q_{-i}$.
\end{theorem}

This follows easily from the lemma above so the proof is omitted.

This demonstrates an important fact about continued fractions: 
that the direction of the sequences of pseudo-squares and 
residues can be reversed (i.e. the indices decrease) by making a 
slight change and applying the same recursive mechanism.  
The presence of one point of symmetry allows a proof that 
another point 
of symmetry exists and that a factorization of $N$ may be obtained 
from this symmetry\footnote{This was actually discovered in the 
opposite order.  It was clear that ambiguous forms that met this 
criteria provided a factorization but was later realized that these 
same forms produced symmetry points.  This was first noticed by Gauss 
\cite{Gauss} and first applied by Shanks \cite{Sh4}.}:

\begin{theorem}
\label{Symmetries}
Let $s = \lfloor \frac{\pi}{2} \rfloor$, where $\pi$ is the 
period from Theorem \ref{strongPeriod}.  If $\pi$ is even then 
(a) $Q_{s+i} = Q_{s-i}$, 
(b) $Q_s \neq Q_0$,
(c) $P_s=P_{s-1}$, and (d) $Q_s \mid 2N$,
for all $i\in \zzz$.
If $\pi$ is odd then, for all $i\in \zzz$, 
\begin{itemize}
\item
$Q_{s+i+1} = Q_{s-i}$, and 

\item
either
(a) $\gcd(Q_s,N)$ is a nontrivial factor of $N$, or 
(b) $-1$ is a quadratic residue of $N$.
\end{itemize}

\end{theorem}

The argument for the first statement is in \cite{W},
pages 641-642. For an elementary proof of both statements, 
see \cite{M}.

\subsection{Infrastructure distance formula}

For $m<n$, and for $\{x_i\}_{i\in \zzz}$, the terms in the 
continued fraction
in (\ref{basicStep}), Shanks defined 
\textbf{infrastructure distance} by

\begin{equation}
\label{distance}
D(x_m,x_n) = \log\left(\prod\limits_{k=m+1}^n x_k\right) .
\end{equation}
We abuse notation and write $D(F_m,F_n)$ as well for this quantity,
where a form $F$ corresponds to a term $x$ in the continued fraction
via the map $x\longmapsto f_+$ (\ref{eqn:CF2QF}).
Lenstra \cite{Len} adds a term of 
$\frac{1}{2}\log(Q_n/Q_m)$ to this (where $Q$ denotes the 
pseudo-square term of $x$), with 
the effect that the resulting formulas are slightly simplified 
but the proofs are more complicated and less intuitive.  
Definition \ref{distance} is also used by Williams in \cite{W}.

Since the quadratic forms are cyclic, in order for the distance 
between two forms to be measured consistently, it must be 
considered modulo the distance around the principal cycle.

\begin{definition} 
\label{regulator}
{\rm Let $\pi$ be the period of the principal cycle.  
The \textbf{regulator} 
$R$ of the class group is the distance around the principal 
cycle, that is},
$R = D(F_0,F_{\pi})$.
\end{definition}

Therefore, distance must be considered modulo $R$, so that 
$D$ is a map from pairs of forms to the interval 
$[0,R) \subset \mathbb{R}$.  The 
addition of two distances must be reduced modulo $R$ as necessary.

\begin{theorem} (infrastructure distance formula)
\label{thrm:14}
If $F_1 \sim F_k$ are
equivalent forms and $G_1 \sim G_\ell$ are equivalent forms and $D_{\rho,1}$ is
the reduction distance for $F_1 * G_1$ and $D_{\rho,2}$ is the reduction
distance for $F_k * G_\ell$ and $m_1$ and $m_k$ are the factors cancelled in each
respective composition (Algorithm \ref{alg:comp}), then

\[
D(F_1\# G_1,F_k\# G_\ell) = D(F_1,F_k)+D(G_1,G_\ell)
+D_{\rho,2}-D_{\rho,1}+\log(m_2/m_1)
\]
\end{theorem}

\pf
Here is a sketch. (For more details, see Theorem A.5.2 in \cite{M}.)

As each quadratic form is associated with a reduced lattice, 
an analysis of
distance requires a connection between reduced lattices
(see \S 3 of \cite{W} for the definition of 
reduced lattice). {\it We use the
notation of Williams \cite{W} without further mention.}

If ${\mathcal L}$ denotes lattice in 
$\qqq (\sqrt{N})$,
let $L({\mathcal L})$ denote the least positive integer
contained in it.

\begin{lemma} (Lemma A.4.2 of \cite{M})
Let $I$ be a primitive ideal and let ${\mathcal L}$ denote the lattice
corresponding to $I$.  If
${\mathcal L}'$ is a lattice with basis $\{1,\xi\}$ and for some 
$\theta$, $\theta
{\mathcal L}' = {\mathcal L}$, then the ideal $J$ corresponding to the 
lattice ${\mathcal L}'$ is a primitive ideal and

\begin{equation}
(L(I)\theta)J = (L(J))I
\end{equation}

\end{lemma}

The method of Voronoi (see for example \cite{W}) is used to obtain a
sequence of adjacent minima, corresponding to a sequence of reduced
lattices.  Consider a sequence of lattices ${\mathcal L}_1$,
${\mathcal L}_2,\cdots$ corresponding to ideals 
$K_1, K_2, \cdots$ corresponding
to binary quadratic forms $F_1,F_2,\cdots$, corresponding to terms
$x_1, x_2, \cdots$ in a continued fraction expansion (\ref{basicStep}). 
If, for two adjacent
lattices in the sequence, $\xi_i$ is defined by 
${\mathcal L}_{i+1} = 1/\xi_i{\mathcal L}_i$, 
then the chain of adjacent minima of ${\mathcal L}_1$ are defined
by $\theta_k = \prod_{i=1}^{k-1}\xi_i$, so $\theta_k{\mathcal L}_k =
{\mathcal L}_1$ (see \cite{W}, \S 3).  Distance between such
lattices is then defined by

\begin{equation}
\label{eq:dist-def}
D({\mathcal L}_k,{\mathcal L}_\ell) = \log(\theta_k/\theta_\ell)
\end{equation}
and this definition of distance corresponds 
exactly to the definition given
for quadratic forms (see \cite{W}, \S 6).

Although this definition has so far only been applied to reduced
ideals (for the definition of reduced ideal, see for example
\cite{W} \S 2) and lattices, the reduction of ideals and lattices
corresponding to quadratic form and continued fraction reduction is
well known:

\begin{lemma} (Lemma A.5.1 in \cite{M}) Let $I$ be any primitive ideal in
$\mathbb{Z}[\sqrt{N}]$.  There exists a reduced ideal $I_k$ and a 
$\theta_k \in I$ such that $(L(I)\theta_k)I_n = (L(I_k))I$.
\end{lemma}

Here $\theta_k$ may be efficiently computed by Voronoi's
method or by continued fractions.  Then the reduction distance is defined by
$D_{\rho} = -\log(\theta_k)$ and may be considered as the distance from $I$ to
$I_k$.

Let $I_1$ denote the ideal corresponding to the 
form $F_1$ in the usual way (as in \cite{Len}),
let $J_1$ be the ideal corresponding to $G_1$, and let
$K_1$ denote the ideal corresponding to $F_1*G_1$.
We have that $(s)K_1 = I_1J_1$, for some $s$.  Let
$K_j$ be a reduced ideal and $\lambda \in K_1$ such that

\begin{equation}
\lambda K_j = K_1.
\end{equation}
Then $K_j$ is the ideal corresponding to $F_1\# G_1$.

Similarly, let $I_k$ denote the ideal corresponding to the
quadratic form $F_k$ and $J_\ell$ be the ideal 
corresponding to the form $G_\ell$.  
If $H_1$ denotes the ideal corresponding to the 
composition $F_k*G_\ell$, then
$(t)H_1 = I_kJ_\ell$, for some $t$.  
Let $H$ be a reduced ideal and choose
$\eta \in H_1$ such that $\eta H = H_1$.
Then $H$ corresponds to $F_k\# G_\ell$.

Let $\mu$ and $\phi$ be such that
$\mu I_k = I_1$ and $\phi J_\ell = J_1$.
Combining these equations, gives

\[
K_j = K_1/\lambda = I_1J_1/\lambda s = (\frac{\mu\phi}{\lambda s})I_kJ_\ell =
(\frac{s\mu\phi}{\lambda t})H_1 = (\frac{s\mu\phi\eta}{\lambda t})H.
\]
Set $\psi = \frac{s\mu\phi\eta}{\lambda t}$ and then $\psi H = K_j$,
so that by (\ref{eq:dist-def}),

\[
D(K_j,H) = -\log(\psi) =
-\log(\mu)-\log(\phi)-\log(\eta)+log(\lambda)-\log(s/t)
\]\[
= D(I_1,I_k)+D(J_1,J_\ell)+D(H_1,H_j)-D(K_1,K_j)+\log(t/s),
\]
as desired.
\qed

\begin{remark}
Shanks stated Square Forms Factorization has an expected 
runtime of $O(\sqrt[4]{N})$ (see Gower \cite{Gower} for a 
detailed discussion of this).

We explain a related idea remarked on by H.
Lenstra \cite{Len}, page 148.

The idea is to first compute the regulator $R$. This has
complexity $O (N^{\frac{1}{5}+\epsilon})$,
assuming the Riemann hypothesis \cite{Len}.
Now use the ``baby-step giant-step'' method (as discussed
in \S 13 of \cite{Len}) to get close to the symmetry point:
\begin{algorithm}
{\small{
{\fontfamily{phv}\selectfont
(Baby-step giant-step)\\
Input: $N$ and $R$\\
Output: Factorization of $N$\\
\begin{enumerate}
\item Compute the form $F$ associated to the first or second
steps of the continued fraction algorithm of
the normalized square root of $N$, (\ref{eqn:normsqrt}).

\item\label{bsgs:1} \vspace{-3mm}\textbf{while} F is not within R/4 of the symmetry point (where
distance is judged using the distance formula in Theorem
\ref{thrm:14}).
\begin{enumerate}
\item \vspace{-2mm}Store $F$ in a Collection $F_c$
\item \vspace{-2mm}$F = F\#F$ (These are the ``giant-steps'')
\end{enumerate}
\item \vspace{-3mm}\label{bsgs:2} Use the intermediate forms in $F_c$ to compose with $F$ until
within $\log N$ of the symmetry point.
\item \vspace{-3mm}\label{bsgs:3}Using the forward and backward steps (see Theorem
\ref{strongPeriod}) of the continued fraction algorithm (``baby
steps''), locate the symmetry point.
\item\vspace{-3mm}
using Lemma \ref{2Symmetries} find a factorization of $N$.
\end{enumerate}
}}}
\end{algorithm}
\end{remark}
Steps \ref{bsgs:1}, \ref{bsgs:2}, and \ref{bsgs:3}, each take $O(\log
N)$, so that the factorization takes $O (N^{\frac{1}{5}+\epsilon})$.

\section{SQUFOF}

Formally, here is the algorithm for factoring $N$:

\begin{algorithm}
\label{alg:fact}
{\small{
{\fontfamily{phv}\selectfont
(SQUFOF)\\
Input: $N$.\\
Output: A factor of $N$\\
\begin{enumerate}
\item \vspace{-3mm}$Q_0 \leftarrow 1, P_0 \leftarrow \lfloor \sqrt{N} \rfloor, Q_1 \leftarrow N-P_0^2$ 
\item \vspace{-3mm}$r \leftarrow \lfloor \sqrt{N} \rfloor$ 
\item \vspace{-3mm}\textbf{while} $Q_i \neq$ perfect square for some $i$ even
\begin{enumerate}
\item \vspace{-2mm}$b_i \leftarrow \left \lfloor \frac{r+P_{i-1}}{Q_i} \right \rfloor$ 
\item \vspace{-1mm}$P_i \leftarrow b_iQ_i-P_{i-1}$ 
\item \vspace{-2mm}$Q_{i+1} \leftarrow Q_{i-1}+b_i(P_{i-1}-P_i)$ 
\item \vspace{-2mm}\textbf{if} $i = 2^n$ for some $n$ Store $(Q_i,2\cdot P_i)$
\end{enumerate}
\item \vspace{-3mm}$F_0 = (\sqrt{Q_i},2\cdot P_{i-1},\frac{P_{i-1}^2-N}{Q_i})$ 
\item \vspace{-3mm}Compose $F_0$ with stored forms according to the binary 
representation of $i/2$ and store result to $F_0$.
\item \vspace{-3mm}$F_0 = (A,B,C)$ 
\item \vspace{-3mm}$Q_0 \leftarrow |A|, P_0 \leftarrow B/2, Q_1 \leftarrow |C|$ 
\item \vspace{-3mm}$q_0 \leftarrow Q_1, p_0 \leftarrow P_0, q_1 \leftarrow Q_0$ 
\item \vspace{-3mm}\textbf{while}\ $P_i \neq P_{i-1}$ and $p_i \neq p_{i-1}$
\begin{enumerate}
\item \vspace{-2mm}Apply same recursive formulas to $(Q_0,P_0,Q_1)$ and $(q_0,p_0,q_1)$
\end{enumerate}
\item \vspace{-3mm}If $P_i = P_{i-1}$, either $Q_i$ or $Q_i/2$ is a nontrivial factor of $N$. 
\item \vspace{-3mm}If $p_i = p_{i-1}$, either $q_i$ or $q_i/2$ is a nontrivial factor of $N$. 
\end{enumerate}}}}
\end{algorithm}

\subsection{Proof}
\label{SummaryProof}

Let $N$, the number to be factored, not be a perfect square.
Expanding the continued fraction for $\sqrt{N}$, let $Q$ be the first
square pseudo-square found on an even index.  Let $r= \sqrt{Q}$.  Let
$F = (r^2,b,c)$ be the associated quadratic form.  Then $(r,b,rc)$,
which reduces with reduction distance $D_\rho = 0$ to $G = (r,b',c')$
is a reduced quadratic form whose square is $F$.  Therefore, by
Theorem \ref{Order2}, $G$ is ambiguous and thus has a symmetry point
in its cycle.

Since by Theorem \ref{thrm:14}, $2D(G_s,G) = D(F_s,F) \pmod R$ where
$F_s$ is the symmetry point of the principal cycle with coefficient
$1$, $D(G_s,G) = D(F_s,F)/2 \pmod {R/2}$.  Since the two points of
symmetry are $R/2$ away from each other, this means that there is a
symmetry point at distance $D(F_s,F)/2$ behind $G$.  Therefore, a
point of symmetry may be found by reversing $G$ and traveling this
short distance.  Now if the coefficient at this symmetry point is $\pm
1$, then there would have been a pseudo-square in the continued
fraction expansion equal to $r$ somewhere before $F$.  If the
coefficient is 2, then this symmetry point could be composed with $G$
to find $2r$ at an earlier point in the principal cycle.  Therefore,
if neither $r$ nor $2r$ were encountered before $F$ in the continued
fraction expansion, then the symmetry point provides a nontrivial
factor for $N$.

\section{Parallel SQUFOF}

With the large amount of computation required for factorization, the
efficiency of a parallel implementation is especially important for
factorization algorithms (see Brent \cite{Br} for a survey and some
terminology).  

There have been proposed two ways to parallelize SQUFOF: using
multipliers and using segments.  We will discuss the segments method
here.  More information on the multipliers method can be found in
Gower \cite{Gower}.

\subsection{Segments}

The segments technique depends upon the ability to use composition to
jump to arbitrary locations in the principal cycle.  The cycle can be
divided into multiple equal-sized sub-sequences and each sub-sequence
can be searched by one of the processors.  As recently as ANTS 2004,
Pomerance suggested investigating parallel SQUFOF (personal
communication; see also \cite{W} page 645).

When factoring using SQUFOF parallelized by segments, we choose a
quadratic form $G$ several steps into the cycle and then square it
several times (how many times is more an art than a science - it
depends on the number of processors and their speed and wanting to
have segments which finish fast but not too fast, say 20-30 in our
case). Call the resulting form $F$. For $i\geq 1$, each $F^{2i}$ is
assigned to processor $i$ as a beginning of another segment,
$[F^{2i}$, $\rho(F^{2i})$, $\rho^2(F^{2i})$, .., $F^{2i+2}]$, where
$\rho$ is the adjacency map.  When processor $i$ finds a pseudo-square
which is a perfect square, that form $H$ may used to find the symmetry
point as follows (Note $H=\rho^{2n}(F^{2i})$, for some $n$).  First,
take the square root of $H$ and reverse it, call this $H'$. This is in
a new cycle of quadratic forms.  Next, compose $H'$ with $F^i$, call
it $H''$. Finally, compose $H''$ with powers of $G$ to
bring it closer to the symmetry point.

\begin{algorithm}
\label{alg:seg-psqufof}
{\small{
{\fontfamily{phv}\selectfont
(Segment-based Parallel SQUFOF)\\
Input: $N$\\
Output: A factor of $N$\\
Preparation:
\begin{enumerate}
\item \vspace{-3mm}$r \leftarrow \lfloor \sqrt{N} \rfloor$ 
\item \vspace{-3mm}$F_0 \leftarrow (1, 2 r, N-r^2)$ 
\item \vspace{-3mm}Cycle $F_0$ several steps forward. 
\item \vspace{-3mm}\textbf{for} $i=1$ to size\ \ \ (size is the logarithmic size of a segment.)
\begin{enumerate}
\item \vspace{-2mm}$F_i \leftarrow F_{i-1}*F_{i-1}$ 
\end{enumerate}
\item \vspace{-3mm}$F \leftarrow F_i$
\end{enumerate}
Processor 0: 
\begin{enumerate}
\item \vspace{-3mm}Assign one processor to search from $F_0$ to $F_{size}$. 
\item \vspace{-3mm}$F_{start} \leftarrow F_{size}$,$F_{end} \leftarrow F_{size}^2$,$F_{rootS} \leftarrow F_{size-1}$,$F_{rootE} \leftarrow F_{size}$,$F_{step} \leftarrow F_{size-1}$
\item \vspace{-3mm}\textbf{while} A factor hasn't been found 
\begin{enumerate}
\item \vspace{-2mm} Wait for a processor to be free and send $F_{start}$, $F_{end}$, and $F_{rootS}$. 
\item \vspace{-2mm}$F_{start} \leftarrow F_{end}$,$F_{rootS} \leftarrow F_{rootE}$,$F_{rootE} \leftarrow F_{rootE}*F_{step}$ ,$F_{end} \leftarrow F_{rootE}^2$ 
\end{enumerate}
\end{enumerate}

\noindent Processor $n$: 
\begin{enumerate}
\item \vspace{-3mm}Receive $F_{start}$, $F_{end}$, and $F_{rootS}$ 
\item \vspace{-3mm}count  $\leftarrow 0$ 
\item \vspace{-3mm}$F_0 = (A,B,C)$ 
\item \vspace{-3mm}\textbf{while} A factor is not found and $F_{start} \ne F_{end}$ 
 \begin{enumerate}
 \item  \vspace{-2mm}Cycle $F_{start}$ forward 2 steps. 
 \item  \vspace{-2mm}count $\leftarrow$ count+1 
 \item  \vspace{-2mm}\textbf{if} $A$ is a perfect square 
  \begin{enumerate}
  \item  \vspace{-1mm}$F_{test} \leftarrow F_{start}^{-1/2}$ 
  \item  \vspace{-1mm}$F_{test} \leftarrow F_{test}*F_{rootS}$ 
  \item  \vspace{-1mm}\textbf{for} $j = size$ to 1 \hspace{0.1in}  (This loop composes $F_{test}$ with the necessary  
   \begin{enumerate}
   \item \vspace{-1mm}\textbf{if} count $> 2^j$ \hspace{0.1in} forms to bring it close to the symmetry point.)
    \item \vspace{-1mm} $F_{test} \leftarrow F_{test}*F_j$ 
    \item  \vspace{-1mm}count $\leftarrow$ count $- 2^j$
   \item  \vspace{-1mm}Search in both directions from $F_{test}$ for a symmetry point. 
   \item  \vspace{-1mm}\textbf{if} Factorization found at symmetry point, output and quit. 
   \end{enumerate}
  \end{enumerate}
 \end{enumerate}
\item \vspace{-3mm}\textbf{if} A factor is still not found, receive new $F_{start}$, 
$F_{end}$, and $F_{rootS}$ and start over.
\end{enumerate}
}}}
\end{algorithm}

Since there is no overlap between the segments searched by the
processors and since the perfect squares appear to be distributed
evenly throughout the principal cycles, this parallelization should be
efficient for any number of processors.  There are two hazards when
choosing selecting the size of the segment.  If the segment size is
too small, the processors will finish their segments so quickly that
receiving new segments will become a bottleneck.  Alternately, if the
segments are too long, the processors may divide up more than the
entire cycle, so that there is overlap.  However, except for rare
numbers that will factor fast regardless, there is significant room in
between these two bounds.
\begin{remark}
The segments based parallelization described here has been implemented
in C using MPI and run on a 64 processor SGI Origin 2800.  Detailed
results and comparisons to the multipliers method can be found in
McMath \cite{M}.  Initial results indicate that the segments method
does indeed continue to be efficient when the number of processors is
increased. 
\end{remark}

The parallelization of SQUFOF by segments involves exactly the same formulas
as the parallelization of the continued fraction factoring algorithm. This
was done in 1987 by Williams andWunderlich [8]. Algorithm 5 of the manuscript
is same as Algorithm 4 of [8], although the former is couched in terms of binary
quadratic forms while the latter uses continued fractions. The equivalence
between binary quadratic forms and continued fractions is well known.

In order to evaluate the efficiency of the segment-based
parallelization, we implemented it and compared it empirically to an
implementation of the multiplier-based version.  The test integers
were all products of randomly chosen primes of roughly equal size.  
Primes of size 80 bit, 100 bit and 120 bit were all tested on 20, 
30, 40, and 50 processors.  This allows an analysis of both how each 
algorithm is affected by the size of the integers and how 
efficiently each algorithm uses an increasing number of processors.


\subsection{Multipliers}

 In 1982, D. Shanks and H. Cohen attempted a
parallelization by having multiple processors attempt to factor $N$,
$3N$, $5N$, etc.  Gower's recent Ph.D. thesis (under S. Wagstaff)
\cite{Gower} analyzed the use of multipliers and found them to be
effective in general but didn't provide much evidence on their
efficiency for parallelization.

The multipliers technique of Gowers-Wagstaff involves generating
multiple version of the factorization algorithm by multiplying $N$ by
products of small square-free numbers $k_i$.  Each product yields a
new number $M_i$ which can be factored on a single processor of a
parallel machine.  If processor $i$ discovers a factor of $M_i$ that
is not from $k_i$, then a factor of $N$ has been found.  
Parallel SQUFOF
using multipliers was considered by Shanks and H. Cohen (when Cohen
visited Shanks at the University of Maryland in 1983, mentioned to the
second author in a private conversation), mentioned by
Williams (\cite{W}, page 645, as an interesting line of research), and
S. Wagstaff and his students (most recently J. Gower \cite{Gower}).

A quick survey of our data for the average runtime shows that for the
segments parallelization, the runtime was cut in half from 20
processors to 50, while the multipliers implementation didn't do quite
so well. The data indicates that the efficient use of multiple
processors for the segments parallelization is roughly unaffected by
increasing the number of processors, while the multipliers
parallelization is less efficient at using a larger number of
processors.  This is the expected result.  As Jason Gower demonstrated
in \cite{Gower}, the use of a multiplier can decrease the runtime by
an average of 27\%.  Therefore, for small numbers of processors, using
multipliers should immediately cut the runtime down.  However, for
larger numbers of processors, the multipliers available aren't used as
efficiently.

Although the data isn't completely clear, the trend is toward segments
being faster than multipliers if enough processors are used.  Based on
the averages, a linear regression predicts a crossover at 80
processors and a quadratic regression predicts a crossover at 47
processors.  The correct answer is probably somewhere within that
range, but even with extensive testing, it would be hard to pin down
the crossover exactly due to the large standard deviations arising in
the data.

\section{Conclusion}
This paper, aside from presenting SQUFOF in its entirety for the first
time, has shown that the algorithm can be presented in terms of an
elegent theoretical framework using two-sided continued fractions and
class groups of quadratic forms over a real quadratic field.  It
further proved the infrastructure distance formula on the cycle of
forms in the class group.  

\subsection*{Acknowledgements}
Daniel Shanks's hand-written notes were kindly made available to the
authors by the executors of his papers (W. Adams, D. Buell, and
H. Williams), to whom we are very grateful.  We are also very grateful
to S. Wagstaff and J. Gowers, who kindly sent us Gower's recent PhD
thesis \cite{Gower}, and to Buell and Williams for many helpful emails.

\end{document}